\documentclass{article}

\usepackage{amsfonts}
\usepackage{latexsym}
\usepackage{amssymb}

\newtheorem{thm}{Theorem}

\newtheorem{prop}{Proposition}
\newtheorem{cor}{Corollary}

\newenvironment{proof}{{\bf Proof.}}{\hfill$\Box$\\[1mm]}

\newcommand{\A}{\mathcal{A}}
\newcommand{\B}{\mathcal{B}}

\renewcommand{\O}{\mathcal{O}}

\newcommand{\U}{\mathcal{U}}

\newcommand{\bbowtie}{\blacktriangleright\!\!\blacktriangleleft}

\begin{document}
\begin{center} 
{\LARGE FRT-duals as Quantum Enveloping Algebras}\\[1cm]
{\large Ulrich Kr\"ahmer}\\[3mm]
Fakult\"at f\"ur Mathematik und Informatik\\
Universit\"at Leipzig, Augustusplatz 10, 04109 Leipzig, Germany\\
E-mail: kraehmer@mathematik.uni-leipzig.de
\end{center}
\abstract{The Hopf algebra generated by the l-functionals on the 
quantum double $\mathbb{C}_q[G] \bowtie \mathbb{C}_q[G]$ is 
considered, where $\mathbb{C}_q[G]$ is the coordinate algebra of a standard 
quantum group and $q$ is not a root of unity. It is shown to be isomorphic
to $\mathbb{C}_q[G]^\mathrm{op} \bowtie U_q(\mathfrak{g})$.
This was conjectured by T. Hodges in [Ho]. As an algebra it can be 
embedded into $U_q(\mathfrak{g}) \otimes U_q(\mathfrak{g})$. Here 
it is proven that
there is no bialgebra structure on
$U_q(\mathfrak{g}) \otimes U_q(\mathfrak{g})$, for which this
embedding becomes a homomorphism of bialgebras. In particular, it is
not an isomorphism.\\
As a preliminary a lemma of [Ho] concerning the
structure of l-functionals on $\mathbb{C}_q[G]$ is generalized. 
For the classical groups a certain choice of root vectors is
expressed in terms of l-functionals. A formula for their coproduct is derived.}

\section{Overview}
\indent Let $\A$ be a coquasitriangular Hopf algebra with universal r-form
$r$ and let $\U(\A)$ be the Hopf subalgebra of the Hopf dual $\A^\circ$
generated by the set of all l-functionals $l^+(a) := r(\cdot \otimes a)$,
$l^-(a) := \bar r(a \otimes \cdot)$, $a \in \A$. We call it the
FRT-dual of $\A$ as it was suggested in [Ho]. There it was shown 
(the finite-dimensional case is treated already
in [Ma]) that there exists an injective algebra homomorphism
$$
	\iota : \U(\A \bowtie \A) \rightarrow \U(\A) \otimes \U(\A)
$$ 
and a surjective Hopf algebra homomorphism
$$
	\zeta : \A^\mathrm{op} \bowtie \U(\A) \rightarrow \U(\A \bowtie \A).
$$
Here $\A \bowtie \B$ denotes the quantum double of the skew-paired 
Hopf algebras $\A$ and $\B$. The skew-pairing of $\A$ and $\A$ in 
$\A \bowtie \A$ is the universal r-form $r$ and the skew-pairing of 
$\U(\A)$ and $\A^\mathrm{op}$ in $\A^\mathrm{op} \bowtie \U(\A)$ is 
the restriction of the canonical pairing of $\A^\circ$ and $\A$. The 
universal r-form on $\A \bowtie \A$ used to define $\U(\A \bowtie \A)$ is 
$\hat r:=\bar r_{41} \bar r_{31} r_{24} r_{23}$ (see section~2 for details).\\
In this paper we continue the investigation of these maps under the 
assumption that $\A$ is the coordinate algebra $\mathbb{C}_q[G]$ of a standard
quantum group associated to a connected complex semi-simple Lie group 
$G$ and $q$ is not a root of unity.\\
The main results are the following facts:
\begin{enumerate}
\item As conjectured in [Ho] $\zeta$ is an isomorphism 
in this case (Theorem~\ref{main}).
\item There exists no bialgebra structure on 
$\U(\mathbb{C}_q[G]) \otimes \U(\mathbb{C}_q[G])$ 
such that $\iota$ becomes a bialgebra 
homomorphism (Theorem~\ref{nogo}).\\
In particular, $\iota$ is not an isomorphism (Corollary~\ref{wuppda}).
\end{enumerate}
We retain the definition of the quantum enveloping algebra 
$U_q(\mathfrak{g})$ from [Ho]. It differs from the usual one
by an extension of the Cartan part. Then $U_q(\mathfrak{g})$ can be 
identified with $\U(\mathbb{C}_q[G])$, if $q$ is not a root of unity
(Proposition~\ref{isom}). Explicitly, one has
$$
	l^+(c^\lambda_{-\mu,\nu}) =
	f^+(c^\lambda_{-\mu,\nu}) K_\mu,\quad 
	l^-(c^\lambda_{-\mu,\nu})=
	f^-(c^\lambda_{-\mu,\nu}) K_{-\nu}
$$ 
for some $f^\pm(c^\lambda_{-\mu,\nu}) \in U_q(\mathfrak{n}_\pm)$ 
(Proposition~\ref{lstru}). Here $c^\lambda_{-\mu,\nu} \in \mathbb{C}_q[G]$ 
are matrix coefficients of the irreducible representation of
$U_q(\mathfrak{g})$ with highest weight $\lambda$.\\  
If $G$ is one of the classical Lie groups,
then there is a choice of the longest word of the
Weyl group of $\mathfrak{g}$, such that the corresponding root vectors of
$U_q(\mathfrak{g})$ occur in the above formula as 
$f^\pm(c^\lambda_{-\mu,\nu})$ for some $c^\lambda_{-\mu,\nu}$
(Proposition~\ref{rootvect}). 
As a corollary one obtains a formula for their coproduct
(Corollary~\ref{copr}).\\[2mm]
There are at least two interpretations of the algebra
$\mathbb{C}_q[G] \bowtie \mathbb{C}_q[G]$.\\
For arbitrary $q$ it is a nonstandard deformation of 
$\mathbb{C}[G \times G]$. In [Ho] it is therefore denoted by 
$\mathbb{C}_q[D(G)]$ where $D(G)$ stands for the double group 
$G \times G$.\\ 
If $q$ is real, it becomes a Hopf $\ast$-algebra which is a
deformation of the algebra of polynomial functions in holomorphic  
and antiholomorphic coordinates on $G$ and describes $G$ as a real Lie
group. It appeared first in this role in the $q$-deformation of the action of 
$SL(2)$ on Minkowski space, see [PW] or [CSSW].\\
Many authors proposed definitions of a quantum enveloping algebra
corresponding to $\mathbb{C}_q[G] \bowtie \mathbb{C}_q[G]$, in
particular, of a $q$-Lorentz algebra. All are based on the requirement that
it should be a Hopf algebra dually paired with 
$\mathbb{C}_q[G] \bowtie \mathbb{C}_q[G]$.\\
One direct approach to such a Hopf algebra is a dualization of the structure of
$\mathbb{C}_q[G] \bowtie \mathbb{C}_q[G]$ in form of a quantum  
codouble $U_q(\mathfrak{g}) \bbowtie U_q(\mathfrak{g})$. It is shown
in [Ma] that $\iota$ would be a Hopf algebra homomorphism into such a 
quantum codouble. Hence it cannot be well-defined by Theorem~\ref{nogo}.\\
In view of the isomorphism $\U(\mathbb{C}_q[G]) \simeq
U_q(\mathfrak{g})$ it seems reasonable to consider the 
FRT-dual $\U(\mathbb{C}_q[G] \bowtie \mathbb{C}_q[G])$ as a
rigorously defined alternative. Theorem~\ref{main} 
is then a dual and purely algebraic form of the Iwasawa decomposition
introduced in [PW] for the
$C^\ast$-completion of $\mathbb{C}_q[G] \bowtie \mathbb{C}_q[G]$.
Note that the images
$\mathrm{Im}\,\hat l^+ \simeq \mathbb{C}_q[G]^\mathrm{op}$ and
$\mathrm{Im}\,\hat l^- \simeq U_q(\mathfrak{g})$
of the l-functionals on
$\mathbb{C}_q[G] \bowtie \mathbb{C}_q[G]$ determine the Iwasawa
decomposition.\\[2mm]
The rest of this paper is divided into three sections: In order
to be self-contained and to fix notations we first recall mainly from [Ho]
and [Ma] some facts about $\A \bowtie \A$ and 
$\U(\A \bowtie \A)$ for an arbitrary coquasitriangular Hopf
algebra $\A$. In the second section we focus on quantum groups and
prove the main results. The last one deals with the relation between 
l-functionals on $\mathbb{C}_q[G]$ and root vectors of
$U_q(\mathfrak{g})$.\\[2mm]
We essentially retain the notations and conventions from [Ho].
We will freely use material that can be found in standard textbooks 
such as [Mo].\\[2mm]
In the original version of this paper only the classical groups were
treated. The author would like to thank T. Hodges and the referee
for pointing out that the proof of Theorem~\ref{main} works with a 
minor modification for arbitrary semi-simple groups. They also
noticed that the well-definedness of 
$U_q(\mathfrak{g}) \bbowtie U_q(\mathfrak{g})$ was an open problem
until now.

\section{Preliminaries on $\A \bowtie \A$ and $\U(\A \bowtie \A)$}
Let $\A$ be a coquasitriangular Hopf algebra with universal r-form $r$.
Then the quantum double $\A \bowtie \A$ is a Hopf algebra
which is the tensor product coalgebra $\A \otimes \A$ endowed with the 
product
$$
	(a \otimes b)(c \otimes d) := 
	(ac_{(2)} \otimes b_{(2)}d)\,
	\bar r(b_{(1)} \otimes c_{(1)})
	r(b_{(3)} \otimes c_{(3)}).
$$
Here $\bar r$ denotes the convolution inverse of $r$ and we use
Sweedlers notation for the coproduct on the right hand side. The
antipode of $\A \bowtie \A$ is given by
$S(a \otimes b):=(1 \otimes S(b))(S(a) \otimes 1)$. See [Ho], [Ma] or [KS1]
for more information about coquasitriangular Hopf algebras and 
quantum doubles.\\
Let $\U(\A)$ be the Hopf subalgebra of the Hopf dual $\A^\circ$ generated by
the set of all l-functionals 
$$
	l^+(a) := r(\cdot \otimes a),\quad l^-(a) := \bar r(a \otimes \cdot),\quad 
	a \in \A.
$$
Following the terminology from [Ho] we call $\U(\A)$ the FRT-dual
of $\A$.\\
If $r$ is a universal r-form on $\A$, then $\bar r_{21}$ is a
universal r-form as well. Note that some formulas in [Ho] differ from
those in this paper1
r because there the latter r-form is used.\\
The Hopf algebra $\A \bowtie \A$ is again coquasitriangular. We define
its FRT-dual $\U(\A \bowtie \A)$ with respect to the universal r-form
$\hat r:=\bar r_{41} \bar r_{31} r_{24} r_{23}$, that is,
\begin{eqnarray}\label{rform}
&& \hat r\left((a \otimes b) \otimes (c \otimes d)\right)\nonumber\\ 
&=&     \bar r(d_{(1)} \otimes a_{(1)}) \bar r(c_{(1)} \otimes a_{(2)})
	r(b_{(1)} \otimes d_{(2)})r(b_{(2)} \otimes c_{(2)})\\ 
&=& \bar r(c_{(1)}d_{(1)} \otimes a)
	r(b \otimes c_{(2)}d_{(2)}).\nonumber
\end{eqnarray}
Consider now the linear maps 
$$
	\theta : \A \bowtie \A \rightarrow \U(\A),\quad
	a \otimes b \mapsto l^+\left(S^{-1}(a)\right)l^-\left(S^{-1}(b)\right)
$$
and
$$
	m : \A \bowtie \A \rightarrow \A,\quad
	a \otimes b \mapsto ab.
$$
Recall that the antipode of a coquasitriangular Hopf algebra is
always bijective, so $\theta$ is well-defined.\\ 
Using the formulas $\bar r(a \otimes b)=r\left(S(a) \otimes b\right)$,
$r(a \otimes b)=r\left(S(a) \otimes S(b)\right)$
(see, e.g. [KS1], Proposition~10.2) and the fact that the coproduct is an algebra
homomorphism one gets
\begin{eqnarray}
	\hat r\left((a \otimes b) \otimes (c \otimes d)\right)
&=& \bar r(c_{(1)}d_{(1)} \otimes a) r(b \otimes c_{(2)}d_{(2)})\nonumber\\
&=& \left\langle l^+\left(S^{-1}(a)\right),(cd)_{(1)}\right\rangle 
	\left\langle l^-\left(S^{-1}(b)\right), (cd)_{(2)}\right\rangle\nonumber\\
&=& \langle \theta(a \otimes b), m(c \otimes d)\rangle.\nonumber 
\end{eqnarray}
For the convolution inverse 
$\bar{\hat r}=\bar r_{23}\bar r_{24}r_{31}r_{41}$ of $\hat r$ one
obtains similarly
$$
	\bar{\hat r}\left((a \otimes b) \otimes (c \otimes d)\right) =
	\left\langle S^{-1} \left(\theta(a \otimes b)\right),
	m(c \otimes d)\right\rangle. 
$$
We denote the l-functionals of $\U(\A \bowtie \A)$ by 
$\hat l^\pm$. The preceding equations imply
$$
	\hat l^+ = \theta^\circ \circ m,\quad
	\hat l^- = m^\circ \circ S^{-1} \circ \theta,
$$
where $\theta^\circ : \A \rightarrow (\A \bowtie \A)^\circ$ and
$m^\circ : \U(\A) \rightarrow (\A \bowtie \A)^\circ$ are linear maps 
dual to $\theta$ and $m$ in the sense
$$
	\langle \theta^\circ(a),b \otimes c \rangle :=
	\langle \theta(b \otimes c), a \rangle,\quad
	\langle m^\circ(f), a \otimes b \rangle := 
	\langle f , ab \rangle.
$$
In particular, the images of $\hat l^+$ and $\hat l^-$ are contained in 
those of $\theta^\circ$ and $m^\circ$, respectively. The map $m$ is
obviously surjective. But 
$S^{-1} \circ \theta : a \otimes b \mapsto l^-(b)l^+(a)$ 
is also surjective by the definition of $\U(\A)$. Hence one even has 
$$
	\mathrm{Im}\,\hat l^+=\mathrm{Im}\,\theta^\circ,\quad
	\mathrm{Im}\,\hat l^-=\mathrm{Im}\,m^\circ.
$$
The definition of $\U(\A \bowtie \A)$ now implies that
the linear map 
$$
	\zeta : \A \otimes \U(\A) \rightarrow 
	\U(\A \bowtie \A), \quad
	\zeta(a \otimes f) := \theta^\circ(a)m^\circ(f) 
$$
is surjective. It is proven in [Ho] that $\zeta$ becomes a Hopf algebra
homomorphism, if one considers $\A \otimes \U(\A)$ with the Hopf
structure $\A^\mathrm{op} \bowtie \U(\A)$.
Here $\A^\mathrm{op}$ denotes the opposite algebra of $\A$ and the
quantum double is constructed with respect to the canonical pairing of
$\U(\A)$ and $\A$.\\
To avoid further notations we will use the product, coproduct and
antipode of $\A$ to express those of $\A^\mathrm{op}$. So the product of 
$a,b \in \A^\mathrm{op}$ is $ba$ and the coproduct and the
antipode of $\A^\mathrm{op}$ are $\Delta$ and $S^{-1}$, respectively.\\
It is also shown in [Ho] that the map 
\begin{equation}\label{embed}
	\iota : \U(\A \bowtie \A) \rightarrow \A^\circ \otimes \A^\circ,\quad
	f \mapsto \langle f_{(1)} , (\cdot \otimes 1) \rangle \otimes 
	\langle f_{(2)} , (1 \otimes \cdot)\rangle
\end{equation}  
is an embedding of algebras and that
$\iota \circ m^\circ = \Delta$ (the coproduct in $\U(\A)$) and 
$\iota \circ \theta^\circ = (l^- \otimes l^+) \circ \Delta$. In particular,  
$\mathrm{Im}\,\iota \subset \U(\A) \otimes \U(\A)$.\\
If $\A$ is finite-dimensional, then any universal r-form $r$ is 
simultaneously a universal R-matrix 
$R$ for the dual Hopf algebra $\A^\circ$ which therefore is quasitriangular.\\
This R-matrix can be used to form a quantum codouble 
$\A^\circ \bbowtie \A^\circ$ of two copies of
$\A^\circ$, see [Ma]. Its structure is completely dual to that of 
$\A \bowtie \A$ - it is the tensor product algebra 
$\A^\circ \otimes \A^\circ$ with a twisted coproduct
\begin{equation}\label{codouble}
	\Delta(a \otimes b) := 
	a_{(1)} \otimes R (b_{(1)} \otimes a_{(2)}) R^{-1} \otimes b_{(2)}.  
\end{equation} 
The map $\iota$ becomes a Hopf algebra homomorphism into
$\A^\circ \bbowtie \A^\circ$.
If $\A$ is in addition factorizable, then both $\iota$ and $\zeta$ are
isomorphisms ([Ma], Theorem~{7.3.5}). As we will see in the next
section, there is no way to
define the above coproduct in a rigorous way for arbitrary
coquasitriangular Hopf algebras $\A$.\\
If $\A$ is a Hopf $\ast$-algebra and $r$ is of real
type, i.e. $r(a^\ast \otimes b^\ast)=\overline{r(b \otimes a)}$, then
$\A \bowtie \A$ is a Hopf $\ast$-algebra with 
involution defined by $(a \otimes b)^\ast:=b^\ast \otimes a^\ast$
([Ma], section~{7.3}, [KS1], section~{10.2.7}).
This applies to the case of the coordinate algebras $\mathbb{C}_q[G]$
treated in the next section if $q$ is real. The involution on $\mathbb{C}_q[G]$ is
the unique one, for which the pairing with the compact real form of 
$U_q(\mathfrak{g})$ ([KS1], section~6.1.7) is a pairing of Hopf $\ast$-algebras.
Then there is a Hopf algebra embedding $a \mapsto 1 \otimes a$ of $\A$ into
$\A \bowtie \A$ and any element of $\A \bowtie \A$ can be written
uniquely as $a^\ast b$ with $a,b \in \A \subset \A \bowtie \A$.
One says that $\A \bowtie \A$ is a realification of $\A$ 
(in [Ma] it is called a complexification).
There is an involution on $\A^\mathrm{op} \bowtie \U(\A)$ defined by
$$
	(a \otimes f)^\ast :=
	(1 \otimes f^\ast)\left(S^2(a)^\ast \otimes 1\right),
$$
for which $\A^\mathrm{op} \bowtie \U(\A)$ becomes a Hopf
$\ast$-algebra and $\zeta$ a $\ast$-homomorphism 
([Ma], Proposition~{7.1.4} and Theorem~{7.3.5}).\\

\section{Application to quantum groups}
We now specialize the preceding considerations to the case where $\A$
is the coordinate algebra of a standard quantum group.\\
Throughout this section $G$ denotes a connected complex semi-simple 
Lie group with Lie algebra $\mathfrak{g}$,
$\{\alpha_i\}_{i=1,\ldots,N}$ a set of simple roots of $\mathfrak{g}$,
$\mathbf{Q}:=\sum_{i=1}^{N} \mathbb{Z} \alpha_i$ the root lattice and
$\mathbf{L}$ the character group of a maximal torus of $G$ which we
identify with a sublattice of the weight lattice of $\mathfrak{g}$
containing $\mathbf{Q}$. For $\lambda,\mu \in \mathbf{L}$ we set 
$\mu < \lambda$ iff $\lambda-\mu$ is a sum of positive roots.
Furthermore, $\langle \cdot , \cdot \rangle$
denotes the scalar product on $\mathbf{L}$ satisfying
$\langle \alpha_i , \alpha_j \rangle = d_i a_{ij}$, where $a_{ij}$ and
$d_i a_{ij}$ are the entries of the Cartan matrix and the symmetrized
Cartan matrix of $\mathfrak{g}$, respectively.\\
We retain the convention from [Ho], where the quantum enveloping
algebra $U_q(\mathfrak{g})$ has generators
$K_\lambda,E_i,F_j$, $\lambda \in \mathbf{L},i,j=1,\ldots,N$
fulfilling the relations 
\begin{eqnarray}
&& K_\lambda K_\mu = K_{\lambda+\mu},\quad
	[E_i,F_j] = \delta_{ij} \frac{K_i-K_i^{-1}}{q^{d_i}-q^{-d_i}},\quad
	K_i:=K_{\alpha_i},\nonumber\\
&& K_\lambda E_i K_\lambda^{-1} =
	q^{\langle \lambda,\alpha_i \rangle} E_i,\quad
	K_\lambda F_j K_\lambda^{-1} = 
	q^{-\langle \lambda,\alpha_j \rangle} F_j\nonumber
\end{eqnarray}
and the $q$-Serre relations ([KS1], eqs. (6.8), (6.9)). The parameter 
$q \in \mathbb{C}\setminus\{0\}$ is assumed to be not a root of
unity. See eqs. (6.19), (6.20) in [KS1] for the definition of the coproduct, the
counit and the antipode of $U_q(\mathfrak{g})$.\\
Note that $U_q(\mathfrak{g})$ as used here is not the most common one,
where $K_\lambda$ is defined only for 
$\lambda \in \mathbf{Q}$. It
depends on the choice of $G$. 
For simply connected $G$ it coincides
with $\check U$ from [Jo] and for the classical groups with
$U_q^\mathrm{ext}(\mathfrak{g})$ from [KS1].\\
There is a $\mathbf{Q}$-grading on $U_q(\mathfrak{g})$ given by
$$
	U_q(\mathfrak{g})=
	\bigoplus_{\lambda \in \mathbf{Q}} U^\lambda_q(\mathfrak{g}),\quad 
	U^\lambda_q(\mathfrak{g}) :=
	\{f \in U_q(\mathfrak{g})\,|\,
	K_\mu f K_{\mu}^{-1} = q^{\langle \lambda,\mu \rangle} f \quad	
	\forall\, \mu \in \mathbf{L}\}.	
$$
Let $U_q(\mathfrak{h}),U_q(\mathfrak{n}_+),U_q(\mathfrak{n}_-)$ 
be the subalgebras generated by the $K_\lambda$, $E_i$ and $F_j$,
respectively. Setting 
$U_q^\lambda(\mathfrak{n}_\pm):=U_q^\lambda(\mathfrak{g}) \cap
U_q(\mathfrak{n}_\pm)$ we have ([Ja], Lemma~{4.12}):
\begin{prop}\label{jantz}
For $f \in U_q^\lambda(\mathfrak{n}_+)$ and 
$g \in U_q^\lambda(\mathfrak{n}_-)$ 
there are $f'_i \in U^{\mu_i}_q(\mathfrak{n}_+)$,
$f''_i \in U^{\lambda-{\mu_i}}_q(\mathfrak{n}_+)$,
$0<\mu_i<\lambda$ and
$g'_j \in U^{\lambda-{\nu_j}}_q(\mathfrak{n}_-)$, 
$g''_j \in U_q^{\nu_j}(\mathfrak{n}_-)$, $0>\nu_j>\lambda$, such that
\begin{eqnarray}
&& \Delta(f) = f \otimes K_\lambda +
	\sum_i f'_i \otimes f''_i K_{\mu_i} +
	1 \otimes f,\nonumber\\ 
&& \Delta(g)= g \otimes 1 +
	\sum_j g'_j K_{\nu_j} \otimes g''_j +
	K_\lambda \otimes g.\nonumber
\end{eqnarray}
\end{prop}
\begin{proof}
We can assume without loss of generality that
$f=E_{i_1} \cdots E_{i_k}$ and $g=F_{i_1} \cdots F_{i_k}$ with 
$\pm \lambda=\alpha_{i_1} + \cdots + \alpha_{i_k}$, 
because $U_q^\lambda(\mathfrak{n}_\pm)$
is spanned by such monomials.
The proof is now an easy induction on $k$. See [Ja] for the details.
\end{proof}
Let $W$ be the Weyl group of $\mathfrak{g}$ generated by the reflections 
$r_i : \alpha_j \mapsto \alpha_j-a_{ij} \alpha_i$. 
Let $E_{\beta_k},F_{\beta_k},k=1,\ldots,n$ be the root vectors of 
$U_q(\mathfrak{g})$ (see [KS1], section~{6.2.3}) associated to the ordering
$\beta_k:=r_{i_1}r_{i_2}\cdots r_{i_{k-1}}\alpha_{i_k}$ of
the set $\mathbf{R}^+$ of positive roots, where 
$r_{i_1}r_{i_2}\cdots r_{i_n}$ is a reduced expression of the longest
element of $W$. Then by the Poincar\'e-Birkhoff-Witt (PBW) theorem the
following monomials form a vector space basis of $U_q(\mathfrak{g})$:
$$
	K_\lambda F_\mathbf{i} E_\mathbf{j} :=
	K_\lambda 
	F_{\beta_1}^{i_1}\cdots F_{\beta_n}^{i_n}
	E_{\beta_1}^{j_1}\cdots E_{\beta_n}^{j_n},\quad
	\lambda \in \mathbf{L},\mathbf{i},\mathbf{j} \in \mathbb{N}_0^n.  
$$
The coordinate algebra $\mathbb{C}_q[G]$ of the standard quantum group
associated to $G$ is the Hopf subalgebra of
$U_q(\mathfrak{g})^\circ$ spanned by the functionals
$c_{u,v}(f):=u(fv)$, $f \in U_q(\mathfrak{g})$, where
$v$ is a vector in the irreducible representation of $U_q(\mathfrak{g})$ 
with highest weight $\lambda \in \mathbf{L}$
and $u$ is a vector in the dual representation, see [Ho]. If 
$\{u_n\},\{v_n\}$ is a pair of dual bases in the representation space
and its dual, then $\Delta(c_{u,v})= \sum_{n} c_{u,v_n} \otimes c_{u_n,v}$. 
If $u,v$ are weight
vectors possessing weights $-\mu,\nu$, then $c_{u,v}$ is denoted by
$c^\lambda_{-\mu,\nu}$ as well.\\
If $G$ is simply connected, $\mathbb{C}_q[G]$ equals $R_q[G]$ from
[Jo]. The relation with $\O(G_q)$ from
[KS1] will be discussed in the next section.\\
The Hopf algebras $\mathbb{C}_q[G]$ are all coquasitriangular. 
A universal r-form $\beta$ is derived in [Ho] from the Rosso form of 
$U_q(\mathfrak{g})$. To be compatible with [FRT] we use 
$r:=\bar\beta_{21}$. This simply exchanges $l^+$ and $l^-$.
It follows from the construction of $r$ that
the l-functionals on $\mathbb{C}_q[G]$ can be identified with elements of
$U_q(\mathfrak{g})$. That is, there is a Hopf algebra embedding of 
$\U(\mathbb{C}_q[G])$ into $U_q(\mathfrak{g})$. 
By Proposition~4.6 in [HLT] this embedding is in fact
surjective. We therefore have:
\begin{prop}\label{isom}
There is an isomorphism 
$\U(\mathbb{C}_q[G]) \simeq U_q(\mathfrak{g})$.
\end{prop}
This was used tacitly in [Ho]. In what follows, we will not distinguish between
$\U(\mathbb{C}_q[G])$ and $U_q(\mathfrak{g})$ any more.\\
In [Ho] the following description of $\mathrm{Ker}\,l^\pm$ was given:
\begin{equation}\label{kernel}
	c_{u,v} \in \mathrm{Ker}\,l^\pm 
	\quad \Leftrightarrow \quad
	u(U_q(\mathfrak{b}_\mp)v)=0.
\end{equation} 
We will use it to prove the next proposition. It 
generalizes Lemma~{3.3} in [Ho].
\begin{prop}\label{lstru}
For $c^\lambda_{-\mu,\nu} \in \mathbb{C}_q[G]$ 
there are $f^\pm(c^\lambda_{-\mu,\nu}) \in U^{\nu-\mu}_q(\mathfrak{n}_\pm)$
with
$$
	l^+(c^\lambda_{-\mu,\nu})= f^+(c^\lambda_{-\mu,\nu})
	K_\mu,\quad
	l^-(c^\lambda_{-\mu,\nu})= f^-(c^\lambda_{-\mu,\nu}) K_{-\nu}. 
$$
\end{prop}
\begin{proof}
We treat only $l^+$, the other case is analogous. Let 
$c^\lambda_{-\mu,\nu}=c_{u,v}$ be given. Fix dual
bases $\{u_n\},\{v_n\}$ as above consisting of weight vectors with weights 
$-\nu_n,\nu_n$, such that $v$ is one of the $v_n$. Let $v'$ be a
highest weight vector and $c^\lambda_{-\mu,\lambda}=c_{u,v'}$. 
Since $l^+$ is a coalgebra homomorphism, we have
\begin{equation}\label{ndoch1}
	\Delta(l^+(c^\lambda_{-\mu,\lambda})) =
	\sum_{n} l^+(c^\lambda_{-\mu,\nu_n}) \otimes 
	l^+(c^\lambda_{-\nu_n,\lambda}).
\end{equation} 
It is known that the proposition holds for $\nu=\lambda$ ([Ho], Lemma 3.3), so
\begin{equation}\label{saegt}
	l^+(c^\lambda_{-\mu,\lambda}) =
	f^+(c^\lambda_{-\mu,\lambda}) K_\mu,\quad
	l^+(c^\lambda_{-\nu_n,\lambda}) =
	f^+(c^\lambda_{-\nu_n,\lambda}) K_{\nu_n}.
\end{equation}
By the first equality and Proposition~\ref{jantz} we can express 
$\Delta(l^+(c^\lambda_{-\mu,\lambda}))$ also as 
\begin{equation}\label{ndoch2}
	f^+(c^\lambda_{-\mu,\lambda}) K_\mu \otimes K_\lambda 
	+ \sum_i 
	f'_i K_\mu \otimes f''_i K_{\xi_i+\mu}
	+ K_\mu \otimes f^+(c^\lambda_{-\mu,\lambda}) K_\mu	
\end{equation} 
with $f'_i \in U_q^{\xi_i}(\mathfrak{n}_+)$,
$f''_i \in U_q^{\lambda-\mu-\xi_i}(\mathfrak{n}_+)$,
$0<\xi_i<\lambda-\mu$. If one compares the $U_q(\mathfrak{h})$-parts
of the terms in 
(\ref{ndoch1}) and (\ref{ndoch2}) in the second tensor component, 
one gets by the second
equality in (\ref{saegt}) and the PBW theorem
$$
	\sum_j l^+(c^\lambda_{-\mu,\nu_{n_j}}) \otimes 
	l^+(c^\lambda_{-\nu_{n_j},\lambda}) =
	\sum_k f'_{i_k} K_\mu \otimes f''_{i_k} K_{\xi_{i_k}+\mu},
$$
where the indices $n_j$ and $i_k$ are those with 
$\nu_{n_i}=\xi_{i_k}+\mu=\nu$.\\
We claim that the elements $l^+(c^\lambda_{-\nu_n,\lambda})$ are
linearly independent. Indeed, assume that there are 
$x_n \in \mathbb{C}$ with
$\sum_n x_n l^+(c_{u_n,v'})=l^+(c_{\sum_n x_n u_n,v'})=0$. 
Since $v'$ is a highest weight vector, (\ref{kernel}) implies
$\sum_n x_n u_n=0$. Hence $x_n=0$ for all $n$, because
$\{u_n\}$ is a basis. 
It follows that all $l^+(c^\lambda_{-\mu,\nu_{n_j}})$ are linear combinations of
$f'_{i_k} K_\mu$. The considered 
$l^+(c^\lambda_{-\mu,\nu})$ is one of them, so the proposition follows.
\end{proof}
Now we are ready to prove the main theorem. 
\begin{thm}\label{main}
There is an isomorphism of Hopf algebras
$$
	\U(\mathbb{C}_q[G] \bowtie \mathbb{C}_q[G]) \simeq 
	\mathbb{C}_q[G]^\mathrm{op} \bowtie U_q(\mathfrak{g}).
$$
\end{thm}
\begin{proof}
It suffices to prove the injectivity of the epimorphism $\zeta$ 
described in section~2. We prove that 
$\zeta':=\iota \circ \zeta : 
\mathbb{C}_q[G]^\mathrm{op} \bowtie U_q(\mathfrak{g}) \rightarrow 
U_q(\mathfrak{g}) \otimes U_q(\mathfrak{g})$ with
$\iota$ from (\ref{embed}) is injective.\\ 
Suppose $f \in \mathrm{Ker}\,\zeta'$, 
$f=\sum_{\lambda \in \mathbf{L},\mathbf{i},\mathbf{j} \in \mathbb{N}_0^n}
a_{\lambda\mathbf{i}\mathbf{j}} \otimes K_\lambda F_\mathbf{i}E_\mathbf{j}$ 
with $a_{\lambda\mathbf{i}\mathbf{j}}=0$ for almost all
$\lambda\mathbf{i}\mathbf{j}$. We have to show that $f$ vanishes.\\
Order $\mathbb{N}_0^n$ in such a way that the weights 
$\mu_\mathbf{j}$ of
$E_\mathbf{j}$ form a nondecreasing (with respect to $<$) sequence.
Let $\mathbf{j}_0$ be the maximal $\mathbf{j}$
for which there exists an $a_{\lambda\mathbf{i}\mathbf{j}} \neq 0$.\\
Recall that $\iota \circ m^\circ = \Delta$ and 
$\iota \circ\theta^\circ = (l^- \otimes l^+) \circ \Delta$. Set
$U_q(\mathfrak{b}_\pm):=U_q(\mathfrak{h})U_q(\mathfrak{n}_\pm)$ and
note that by Proposition~\ref{lstru} and Proposition~\ref{jantz} we have
$$
	(l^- \otimes l^+) \circ \Delta(a_{\lambda \mathbf{i} \mathbf{j}}) \in 
	U_q(\mathfrak{b}_-) \otimes U_q(\mathfrak{b}_+),\quad
	\Delta(K_\lambda F_\mathbf{i}) \in 
	U_q(\mathfrak{b}_-) \otimes U_q(\mathfrak{b}_-).
$$
Hence only $\Delta(E_\mathbf{j})$ contribute to the 
$U_q(\mathfrak{n}_+)$-part in the first tensor component.
Expand them according to Proposition~\ref{jantz}. Then the
PBW theorem implies that
$$
	\sum_{\lambda\mathbf{i}}
	(l^- \otimes l^+) \circ \Delta(a_{\lambda\mathbf{i}\mathbf{j}_0}) \cdot
	\Delta(K_\lambda F_\mathbf{i}) \cdot 
	(E_{\mathbf{j}_0} \otimes K_{\mu_{\mathbf{j}_0}})
$$ 
is linearly independent from the other terms occuring in $\zeta'(f)$ 
and vanishes separately. Since 
$U_q(\mathfrak{g}) \otimes U_q(\mathfrak{g})=
U_q(\mathfrak{g} \oplus \mathfrak{g})$ is free of zero divisors
([DK], Corollary~{1.8}), we get
$$
	\sum_{\lambda \mathbf{i}}
	(l^- \otimes l^+) \circ \Delta(a_{\lambda \mathbf{i} \mathbf{j}_0}) \cdot
	\Delta(K_\lambda F_\mathbf{i}) = 0.
$$
The same argument applied to the maximal $\mathbf{i}_0$
and the second tensor component shows
$$
	\sum_\lambda
	(l^- \otimes l^+) \circ \Delta(a_{\lambda \mathbf{i}_0 \mathbf{j}_0}) \cdot
	\Delta(K_\lambda) = 0.
$$
By Proposition~\ref{lstru} we can write 
$(l^- \otimes l^+) \circ \Delta(a_{\lambda \mathbf{i}_0 \mathbf{j}_0})$ as
$\sum_{\xi \in \mathbf{L}} u_{\lambda\xi}K_\xi \otimes v_{\lambda\xi}K_{-\xi}$ 	 
for some 
$u_{\lambda\xi} \otimes v_{\lambda\xi} \in U_q(\mathfrak{n}_-) \otimes 
U_q(\mathfrak{n}_+)$. Then the last equation becomes
$$
	\sum_{\xi\lambda}
	u_{\lambda\xi}K_{\xi+\lambda} \otimes v_{\lambda\xi}K_{-\xi+\lambda} = 0.
$$
This implies $u_{\lambda\xi} \otimes v_{\lambda\xi} = 0$ for all
$\lambda,\xi$. Finally,
$$
	(l^- \otimes l^+) \circ \Delta(a_{\lambda \mathbf{i_0j_0}})=0
$$
implies $a_{\lambda \mathbf{i_0j_0}}=0$ in contradiction with the
assumption, because $(l^- \otimes l^+) \circ \Delta$ is injective by the
definition of $\U(\mathbb{C}_q[G])$.
\end{proof}
In contrast to their h-adic counterparts $U_h(\mathfrak{g})$ 
defined over the ring of formal power series $\mathbb{C}[[h]]$, the
Hopf algebras $U_q(\mathfrak{g})$ over $\mathbb{C}$ are not 
quasitriangular. Nevertheless, parts of the theory of
$U_h(\mathfrak{g})$ carry over to $U_q(\mathfrak{g})$, since 
the l-functionals encode the R-matrix of $U_h(\mathfrak{g})$
to some extent. Hence it is not a priori clear that there is no way
to define the twisted coproduct (\ref{codouble}) as well on 
$U_q(\mathfrak{g}) \otimes U_q(\mathfrak{g})$. But we show now 
that this is in fact impossible.
\begin{thm}\label{nogo}
There exists no bialgebra structure on 
$U_q(\mathfrak{g}) \otimes U_q(\mathfrak{g})$ such that $\iota$
becomes a homomorphism of bialgebras. 
\end{thm}
\begin{proof}
Suppose that the opposite holds. Then 
$\iota \circ \theta^\circ$ is a bialgebra homomorphism as well. Note that
$\pm(\nu-\mu) \notin \sum_{i=1}^N \mathbb{N}_0 \alpha_i$ implies
$l^\pm(c^\lambda_{-\mu,\nu})=0$ by (\ref{kernel}). Using this and
Proposition~\ref{lstru} one computes 
\begin{eqnarray}
	\Delta(K_{\lambda} \otimes K_{-\lambda})
&=& \Delta \circ \iota \circ
	\theta^\circ(c^\lambda_{\lambda,-\lambda})\nonumber\\
&=& (\iota \circ \theta^\circ \otimes \iota \circ \theta^\circ) \circ 
	\Delta(c^\lambda_{\lambda,-\lambda})\nonumber\\
&=& \sum_{n} K_\lambda \otimes 
	f^+(c^\lambda_{\lambda,\nu_n}) K_{-\lambda}
	\otimes f^-(c^\lambda_{-\nu_n,-\lambda}) K_{\lambda} \otimes
	K_{-\lambda}.\nonumber
\end{eqnarray}
This must be an invertible element of $U_q(\mathfrak{g})^{\otimes 4}$,
because $\Delta$ is an algebra homomorphism and 
$K_{\lambda} \otimes K_{-\lambda}$ is invertible. Since 
$K_\lambda \otimes K_{-\lambda} \otimes K_{\lambda} \otimes
K_{-\lambda}$ is invertible, 
$\sum_{n} f^+(c^\lambda_{\lambda,\nu_n}) \otimes 
f^-(c^\lambda_{-\nu_n,-\lambda})$ is an invertible element of 
$U_q(\mathfrak{g})^{\otimes 2}$.\\
An invertible element of a graded algebra must be homogeneous -
the product of the homogeneous components of highest
degrees $n_0,m_0$ of the element and its inverse must be of degree
zero, so $m_0=-n_0$, the same must hold for the components of lowest
degrees $n_1,m_1$, so $m_1=-n_1$ and
$n_1 \le n_0$ and $m_1 \le m_0$ implies then $m_0=m_1=-n_0=-n_1$.
By Proposition~\ref{lstru}
$\sum_{n} f^+(c^\lambda_{\lambda,\nu_n}) \otimes 
f^-(c^\lambda_{-\nu_n,-\lambda})$ is not homogeneous with respect to
the $\mathbf{Q} \times \mathbf{Q}$-grading of 
$U_q(\mathfrak{g}) \otimes U_q(\mathfrak{g})$, so 
we obtain a contradiction.
\end{proof}
\begin{cor}\label{wuppda}
The map $\iota$ is not surjective.
\end{cor}

\section{L-functionals and root vectors}
The root vectors of $U_q(\mathfrak{g})$ are defined in terms of the
action of the braid group of $\mathfrak{g}$ on $U_q(\mathfrak{g})$.
Since this action is not given by coalgebra homomorphisms, 
it is not possible to compute their coproduct directly from their
definition.\\
However, it is mentioned in [KS1] on p. 278 that for $G=SL(N+1)$ there is
a choice of $r_{i_1}r_{i_2}\cdots r_{i_n}$, such that the root vectors
are certain $f^\pm(c^\lambda_{-\mu,\nu})$ from
Proposition~\ref{lstru}. This allows to compute their coproduct explicitly.\\
In this section we generalize this result to the
other classical Lie groups. The main tool will be the following proposition:
\begin{prop}\label{soibel}
For $i<j$ there are $x_{ij}(\mathbf{k}),y_{ij}(\mathbf{k}) \in \mathbb{C}$, 
such that
\begin{eqnarray}
&& E_{\beta_i} E_{\beta_j} - 
	q^{\langle \beta_i , \beta_j \rangle}
	E_{\beta_j} E_{\beta_i}=
	\sum_{\mathbf{k} \in \mathbb{N}_0^{j-i-1}}
	x_{ij}(\mathbf{k}) E_{\beta_{i+1}}^{k_1}\cdots E_{\beta_{j-1}}^{k_{j-i-1}},
	\label{endli}\\
&& F_{\beta_i} F_{\beta_j} - 
	q^{-\langle \beta_i , \beta_j \rangle}
	F_{\beta_j} F_{\beta_i}=
	\sum_{\mathbf{k} \in \mathbb{N}_0^{j-i-1}}
	y_{ij}(\mathbf{k}) F_{\beta_{i+1}}^{k_1}\cdots F_{\beta_{j-1}}^{k_{j-i-1}}.
	\label{endli2}
\end{eqnarray}
If $\beta_i+\beta_j \neq \sum_{l=1}^{j-1} k_l \beta_{i+l}$,
then $x_{ij}(\mathbf{k})=y_{ij}(\mathbf{k})=0$.%
\end{prop}
\begin{proof} 
The two relations (\ref{endli}), (\ref{endli2}) are proven in [KS2], 
Theorem~{3.2.3}.\\ 
Conjugating (\ref{endli}) with $K_\lambda$ one gets
$$
	\sum_{\mathbf{k}}
	(q^{\langle \lambda,\beta_i+\beta_j \rangle}-
	q^{\langle \lambda,k_1\beta_{i+1}+\cdots+k_{j-i-1}\beta_{j-1}\rangle})
	x_{ij}(\mathbf{k}) 
	E_{\beta_{i+1}}^{k_1}\cdots E_{\beta_{j-1}}^{k_{j-i-1}} = 0.
$$
The PBW theorem implies 
$q^{\langle \lambda,\beta_i+\beta_j \rangle}=
q^{\langle \lambda,k_1\beta_{i+1}+\cdots+k_{j-i-1}\beta_{j-1}\rangle}$
or $x_{ij}(\mathbf{k})=0$. Since $q$ is not a root of
unity and $\lambda$ was arbitrary, the additional claim follows for the
$x_{ij}(\mathbf{k})$. The same argument applies to the $y_{ij}(\mathbf{k})$. 
\end{proof}
We will use a special ordering of the positive roots, in which
most if not all terms on the right hand side of (\ref{endli}),
(\ref{endli2}) vanish. To define it, we
first arrange the positive roots in the following way as parts of matrices:
$$
	\beta_{ij} =
	\left\{
	\begin{array}{lll}
	\begin{array}{l}
	\sum_{k=i}^{j-1} \alpha_k 
	\end{array}
	& 
	& \mathfrak{g}=\mathfrak{sl}_{N+1},\\
	\begin{array}{l}
	\sum_{k=i}^{j-1} \alpha_k\\
	\sum_{k=i}^{N} \alpha_k +
	\sum_{k=j'}^{N} \alpha_k
	\end{array}
	&
	\begin{array}{l}
	j \le N+1\\
	j>N+1 \end{array} 
	& \mathfrak{g}=\mathfrak{so}_{2N+1},\\
	\begin{array}{l}
	\sum_{k=i}^{j-1} \alpha_k\\
	\sum_{k=i}^{N} \alpha_k +
	\sum_{k=j'}^{N-1} \alpha_k
	\end{array}
	& 
	\begin{array}{l}
	j \le N+1\\
	j>N+1
	\end{array}
	& \mathfrak{g}=\mathfrak{sp}_{2N},\\
	\begin{array}{l}
	\sum_{k=i}^{j-1} \alpha_k\\
	\sum_{k=i}^{N-2} \alpha_k+\alpha_N\\
	\sum_{k=i}^{N} \alpha_k\\
	\sum_{k=i}^{N} \alpha_k +
	\sum_{k=j'}^{N-2} \alpha_k
	\end{array}
	&
	\begin{array}{l}
	j \le N\\
	j = N+1\\
	j = N+2\\
	j>N+2
	\end{array}
	& \mathfrak{g}=\mathfrak{so}_{2N},
\end{array}\right.  
$$
where 
$j':=2N+2-j$ for $\mathfrak{g}=\mathfrak{so}_{2N+1}$, 
$j':=2N+1-j$ for $\mathfrak{g}=\mathfrak{sp}_{2N},\mathfrak{so}_{2N}$
and the indices take the values
$$
	\begin{array}{ll}
	i=1,\ldots,N,j=i+1,\ldots,N+1 &
	\mathfrak{g}=\mathfrak{sl}_{N+1},\\
	i=1,\ldots,N,j=i+1,\ldots,i'-1 &
	\mathfrak{g}=\mathfrak{so}_{2N+1},\\
	i=1,\ldots,N,j=i+1,\ldots,i' &
	\mathfrak{g}=\mathfrak{sp}_{2N},\\
	i=1,\ldots,N-1,j=i+1,\ldots,i'-1 &
	\mathfrak{g}=\mathfrak{so}_{2N}.
	\end{array}  
$$
Now we fix the expression $\prod_{k=N}^1 a_k$ for the longest word of $W$, where 
$$
	a_k := 
	\left\{\begin{array}{ll}
	\begin{array}{l}
	\prod_{i=1}^k r_i 
	\end{array} &
	\mathfrak{g}=\mathfrak{sl}_{N+1},\\
	\begin{array}{l}
	(\prod_{i=k}^N r_i) (\prod_{j=N-1}^k r_j) 
	\end{array} &
	\mathfrak{g}=\mathfrak{so}_{2N+1},\mathfrak{sp}_{2N},\\
	\begin{array}{ll}
	1 & k=N\\
	(\prod_{i=k}^{N-2} r_i) 
	r_N (\prod_{j=N-1}^k r_j) & N-k \neq 0 \mbox{ odd}\\
	(\prod_{i=k}^{N-1} r_i) r_N (\prod_{j=N}^k r_j) & N-k \neq 0 \mbox{ even}
	\end{array} 
	&
	\mathfrak{g}=\mathfrak{so}_{2N}. 
	\end{array}\right.
$$
Then the induced ordering $\prec$ of $\mathbf{R}^+$ is as follows: 
$$
	\beta_{ij} \prec \beta_{kl} 
	\quad \Leftrightarrow \quad 
	\left\{
	\begin{array}{ll}
	\begin{array}{l}
	i<k \mbox{ or } i=k,j<l 
	\end{array}
	& \mathfrak{g}=\mathfrak{sl}_{N+1},\\
	\begin{array}{l}
	k<i \mbox{ or } i=k,l<j
	\end{array}
	& 
	\mathfrak{g}=\mathfrak{so}_{2N+1},\mathfrak{so}_{2N},\\
	\begin{array}{l}
	k<i \mbox{ or } i=k,j=N+1 \mbox{ or}\\
	i=k,l<j,j \neq N+1,l \neq N+1
	\end{array}  
	& \mathfrak{g}=\mathfrak{sp}_{2N}.
\end{array}\right.  
$$
Originally the quantum group coordinate algebras were defined only
for the classical groups in terms of generators and relations [FRT]. 
The generators are the matrix coefficients $u^i_j$ of the vector representation of
$U_q(\mathfrak{g})$ (the first fundamental representation which
defines $\mathfrak{g}$ as a matrix Lie algebra) with respect to some basis.
For the relations we refer to chapter~9 of [KS1]. There the resulting
Hopf algebras are denoted by $\O(G_q)$.\\
If $q$ is not a root of unity, then  
$\O(G_q)$ defined in this way is isomorphic to 
$\mathbb{C}_q[G]$ as used in the last section
for all $G$ except $G=SO(2N+1)$. In this case one has $\O(G_{q^2}) \simeq
\mathbb{C}_q[G]$.\\
This is a consequence of the 
Peter-Weyl theorem ([KS1], Theorem~{11.22}). 
The latter is stated in [KS1] under
the assumption that $q$ is transcendental. According to Remark~3 on
p. 415 of [KS1] and Corollaries~{4.15} and~{5.22} from [LR] the result holds also 
for $q$ not a root of unity.\\
We abbreviate $f \sim g$ iff $f=x g$ with some 
$x \in \mathbb{C}\setminus\{0\}$ and $(l^\pm)^i_j:=l^\pm(u^i_j)$. 
Then the following statement holds:
\begin{prop}\label{rootvect}
If $ij$ appear as indices of a positive root $\beta_{ij}$, then
$$
	(l^+)^i_j \sim (l^+)^i_i E_{\beta_{ij}},\quad
	(l^-)^j_i \sim (l^-)^i_i F_{\beta_{ij}},
$$
except if $\mathfrak{g}=\mathfrak{sp}_{2N}$ and $j=i'$. 
In this case, there are $x,y \in \mathbb{C}$, such that
\begin{eqnarray}
&& (l^+)^i_{i'} \sim
	(l^+)^i_i (E_{\beta_{ii'}}-x E_{\beta_{ii'-1}}E_i),\nonumber\\
&& (l^-)^{i'}_i \sim
	(l^-)^i_i (F_{\beta_{ii'}}-y F_{\beta_{ii'-1}}F_i).\nonumber
\end{eqnarray}
\end{prop}
\begin{proof}
Since this is known for $\mathfrak{g}=\mathfrak{sl}_{N+1}$, 
we consider only the remaining cases. We also will
consider only the $E_{\beta_{ij}}$. The $F_{\beta_{ji}}$ are
treated similarly.\\
The proof is by induction over $j-i$. 
By the lists of $(l^\pm)^i_j$ in section~{8.5.2} of [KS1] 
the claim holds for $j-i=1$.
All occurring $(l^+)^i_j$ except $(l^+)^{N-1}_{N+1}$ for 
$\mathfrak{g}=\mathfrak{sp}_{2N},\mathfrak{so}_{2N}$ 
can be calculated from the recurrence relation 
\begin{equation}\label{jetze1}
	(q-q^{-1})(l^+)^i_j = 
	- \left[(l^+)^i_k,(l^+)^k_j\right] (l^-)^k_k. 
\end{equation} 
Here $k$ with $i<k<j$ is arbitrary with $k \neq i',j'$
([KS1], Proposition~{8.29}).\\
We choose $k=j-1$. This is admissible in all cases except 
$\mathfrak{g}=\mathfrak{sp}_{2N},\mathfrak{so}_{2N}$ and 
$j=N+1$. These must be treated separately afterwards.\\  
By the explicit lists of the $(l^\pm)^i_j$ in [KS1] there are
$\lambda_k \in \mathbf{L}$ such that
$$
	(l^\pm)^k_k=K_{\pm\lambda_k},\quad 
	(l^+)^{j-1}_j \sim (l^+)^{j-1}_{j-1}  E_{f(j-1)},\quad
	f(k):=  \left\{
	\begin{array}{ll}
	k & k \le N,\\
	k'-1 & k>N.
	\end{array}\right.  
$$ 
Inserting this and the induction hypothesis into (\ref{jetze1}) we get
\begin{equation}\label{brumm}
	(l^+)^i_j \sim
	(l^+)^i_i (E_{f(j-1)}E_{\beta_{ij-1}}-q^{-g(i,j-1)}E_{\beta_{ij-1}}E_{f(j-1)})
\end{equation} 
with $g(i,j-1)=\langle \lambda_{j-1},\beta_{ij-1}\rangle-
\langle \lambda_i,\alpha_{f(j-1)}\rangle$.\\
Inserting the explicit formulas for $\lambda_k,\beta_{ij}$,
$\langle \alpha_i,\alpha_j \rangle$ one gets after some lengthy calculations
\begin{equation}\label{g}
	g(i,j-1)=\left\{
	\begin{array}{ll}
	2\quad & \mathfrak{g}=\mathfrak{sp}_{2N},j=i'\\
	-\langle \alpha_{f(j-1)},\beta_{ij-1}\rangle\quad &
	\mbox{otherwise}
	\end{array}\right..
\end{equation} 
In our ordering of $\mathbf{R}^+$ we have 
$\alpha_{f(j-1)} \prec \beta_{ij-1}$ for $i<j-1$ and 
$j \neq i'$ which holds in all cases except
$\mathfrak{g}=\mathfrak{sp}_{2N},j=i'$. Since
$\mathfrak{g}=\mathfrak{so}_{2N},j=N+1$ was excluded we furthermore 
have $\alpha_{f(j-1)}+\beta_{ij-1}=\beta_{ij}$  
and there is no other linear combination of roots between 
$\alpha_{f(j-1)}$ and $\beta_{ij-1}$ equal to $\beta_{ij}$. Hence the
exponent in (\ref{brumm}) is in all
considered cases except $\mathfrak{g}=\mathfrak{sp}_{2N},j=i'$ 
the same as the one which appears on the left hand side
of (\ref{endli}) and the claim reduces to Proposition~\ref{soibel}
(note that for the classical groups $(l^+)^i_j \neq 0$ for all 
$i \le j$, as follows for example from [KS1], Theorem~{8.33}).\\
For $\mathfrak{g}=\mathfrak{sp}_{2N},j=i'$ we obtain
$$
	(l^+)^i_{i'} \sim
	(l^+)^i_i (E_iE_{\beta_{ii'-1}}-q^{-2}E_{\beta_{ii'-1}}E_i) \sim 
	(l^+)^i_i (E_{\beta_{ii'}}-x E_{\beta_{ii'-1}}E_i)
$$
for some $x \in \mathbb{C}$, 
because $\langle \alpha_i,\beta_{ij-1} \rangle=0$.\\
It remains to treat the excluded cases $(l^+)^i_{N+1}$ for 
$\mathfrak{g}=\mathfrak{sp}_{2N},\mathfrak{so}_{2N}$.\\
By the explicit lists of $(l^+)^i_j$ in [KS1] we have 
for $\mathfrak{g}=\mathfrak{sp}_{2N}$, $i=N-1$ 
\begin{eqnarray}
	(l^+)^{N-1}_{N+1} 
&\sim& (l^+)^{N-1}_{N-1}(E_NE_{N-1}
	-q^{-2}E_{N-1}E_N)\nonumber\\
&=& (l^+)^{N-1}_{N-1}(E_NE_{N-1}
	-q^{\langle \alpha_{N-1},\alpha_N \rangle}E_{N-1}E_N)\nonumber\\
&\sim& (l^+)^{N-1}_{N-1}E_{\alpha_{N-1}+\alpha_N}\nonumber
\end{eqnarray}
by the same argument as above. For $\mathfrak{g}=\mathfrak{so}_{2N}$
the lists directly contain
$$
	(l^+)^{N-1}_{N+1} \sim (l^+)^{N-1}_{N-1}E_N,
$$
so the proposition holds in these cases.\\
For $i < N-1$ we need a second induction on $i$ starting with 
$i=N-1$. We again use the recurrence relation (\ref{jetze1}),
but now with $k=i+1$ (which is possible for $i<N-1$) getting
by induction
\begin{eqnarray}
	(l^+)^i_{N+1}
&\sim&[(l^+)^{i}_{i+1},(l^+)^{i+1}_{N+1}](l^-)^{i+1}_{i+1} \nonumber\\
&\sim&(l^+)^i_i (E_iE_{\beta_{i+1N+1}}-
	q^{\langle \lambda_{i+1},\alpha_i \rangle-
	\langle \lambda_i,\beta_{i+1N+1} \rangle}
	E_{\beta_{i+1N+1}}E_i).\nonumber
\end{eqnarray} 
In all cases $\langle \lambda_{i+1},\alpha_i \rangle=1$ and the second
term in the exponent vanishes, since in $\beta_{i+1N+1}$ only
$\alpha_j$ with $j>i$ occur. Since 
$\langle \alpha_i,\beta_{i+1N+1}\rangle=-1$ and
$\alpha_i \succ \beta_{i+1N+1}$, the same
argumentation as above yields
\begin{eqnarray}
	(l^+)^{i}_{N+1}
&\sim& (l^+)^i_i (E_{\beta_{i+1N+1}}E_i-
	q^{\langle \alpha_i,\beta_{i+1N+1} \rangle}
	E_iE_{\beta_{i+1N+1}})\nonumber\\
&\sim&(l^+)^i_i E_{\beta_{iN+1}}.\nonumber
\end{eqnarray} 
\end{proof}
\begin{cor}\label{copr}
We have
\begin{eqnarray}
&& \Delta(E_{\beta_{ij}}) 
	\sim ((l^-)^i_i \otimes (l^-)^i_i) 
	\sum_{k=i}^{j}  (l^+)^i_k \otimes (l^+)^k_j, \nonumber\\
&& \Delta(F_{\beta_{ij}}) 
	\sim ((l^+)^i_i \otimes (l^+)^i_i) 
	\sum_{k=i}^{j}  (l^-)^j_k \otimes (l^-)^k_i \nonumber
\end{eqnarray}
 except when $\mathfrak{g}=\mathfrak{sp}_{2N}$ and $j=i'$. In this
case, we have
\begin{eqnarray}
&& \Delta(E_{\beta_{ii'}}) \sim
	x \Delta(E_{\beta_{ii'-1}}E_i)+((l^-)^i_i \otimes (l^-)^i_i) 
	\sum_{k=i}^{i'}  (l^+)^i_k \otimes (l^+)^k_{i'},\nonumber\\
&& \Delta(F_{\beta_{ii'}}) \sim
	y \Delta(F_{i\beta_{i'-1}}F_i)+
	((l^+)^i_i \otimes (l^+)^i_i) 
	\sum_{k=i}^{i'} (l^-)^{i'}_k \otimes (l^-)^k_i.\nonumber
\end{eqnarray} 
\end{cor}

\end{document}